\documentclass[A4paper, 12pt]{amsart}

\usepackage{amsmath, amssymb}
\usepackage{amsthm}
\usepackage{multicol}

\newtheorem{thm}{Theorem}%[section]
\newtheorem{prop}[thm]{Proposition}

\newtheorem{cor}[thm]{Corollary}

\theoremstyle{definition}
\newtheorem{defi}[thm]{Definition}
\newtheorem{ex}[thm]{Example}
\newtheorem{as}[thm]{Assumption}

\theoremstyle{remark}
\newtheorem{rem}[thm]{Remark}

\renewcommand{\t}{\tau}
\newcommand{\tilV}{\tilde{V}}
\newcommand{\tilNu}{\mathcal{\widetilde V}}
\newcommand{\vol}{{\rm Vol}}
\newcommand{\dert}[1]{\frac{\partial #1}{\partial t}}
\newcommand{\dertau}[1]{{\partial #1\over\partial\tau}}
\newcommand{\Ric}{{\rm Ric}}
\newcommand{\Hess}{{\rm Hess}\,}
\newcommand{\tr}{{\rm tr}}
\renewcommand{\L}{\mathcal{L}}
\renewcommand{\l}{\ell}
\newcommand{\I}{\mathcal{I}}
\newcommand{\R}{\mathbb{R}}
\renewcommand{\[}{\Bigl[}
\renewcommand{\]}{\Bigr]}
\renewcommand{\(}{\Bigl(}
\renewcommand{\)}{\Bigr)}
\renewcommand{\|}{\Big|}
\renewcommand{\iint}{\int\!\!\!\int}
\renewcommand{\phi}{\varphi}
\renewcommand{\epsilon}{\varepsilon}
\newcommand{\bt}{{\bar{\tau}}}

\title{On the asymptotic reduced volume of the Ricci flow}
\author{Takumi Yokota}
\address{Graduate School of Pure and Applied Sciences\\
University of Tsukuba\\
305-8571 Tsukuba\\
Japan}
\email{takumiy@math.tsukuba.ac.jp}           %  \\
\date{April 6, 2009}
\begin{document}
\maketitle

\begin{abstract}
In this paper,
we consider two different monotone quantities defined for the Ricci flow
and show that their asymptotic limits coincide for any ancient solutions.
One of the quantities we consider here is Perelman's reduced volume,
while the other is the local quantity discovered by Ecker, Knopf, Ni and Topping.
This establishes a relation between these two monotone quantities. 
\end{abstract}

\section{Introduction}
A smooth one-parameter family of Riemannian metrics $(M, g(t)), t \in I$
with $I \subset \R$ being an interval,
is called  a {\it Ricci flow} when it evolves along the equation
\begin{equation}
\dert{}g = -2\Ric(g(t)),
\end{equation}
where $\Ric(g(t))$ denotes the Ricci tensor of $g(t)$.
In this paper,
we are concerned with ancient solutions and
will use $\t := -t$ to represent the reverse time parameter.
An {\it ancient solution} is a Ricci flow $(M, g(t))$ which exists for all $t \in (-\infty, 0]$.

In the study of the Ricci and other geometric flows,
monotone quantities have been playing significant roles.
One of Perelman's achievements in his seminal paper \cite{P} is the monotonicity of the reduced volume given by
$$\tilV_{(p, 0)}(\t) := \int_M (4\pi\t)^{-n/2} e^{-\l(\cdot, \t)} d\mu_{g(\t)} \text{ for } \t>0.$$
Here, and henceforth, $d\mu$ denotes the volume element.
The precise definition of the reduced distance $\l = \l_{(p, 0)}$
from the base point $(p, 0) \in M\times \{0\}$ is recalled in the next section.
An application of the monotonicity of the reduced volume is
to rule out local collapsing for the Ricci flow (\cite[\S7]{P}).
Throughout this paper,
we adopt the convention that $\tilV_{(p, 0)}(\t) \equiv 1$ for (and only for) the {\it Gaussian soliton},
i.e., the trivial Ricci flow $(\R^n, g(\t) \equiv g_E)$ on the Euclidean space.

Subsequently,
another monotone quantity of the form
$$I_{(p, 0)}(r) := \frac{1}{r^n} \iint_{E_r} \[|\nabla \l|^2 + R\(n\log \frac{r}{\sqrt{4\pi\t}} -\l\)\] d\mu d\t \text{ for } r>0$$
was also discovered by Ecker, Knopf, Ni and Topping (\cite{EKNT}).
Here $R=R(\cdot, \t)$ denotes the scalar curvature of $g(\t)$
and $E_r$ is a certain subset, called `pseudo heat ball', of the space-time.
The precise definition is given in the next section as well.

At this point, it is natural to ask the following question (\cite{EKNT}).
\begin{quotation}
{\it How is the local monotone quantity $I(r)$ related to the global one $\tilV(\t)$?}
\end{quotation}

Partially motivated by this question,
we now state the main result of the present paper as follows:
\begin{thm}\label{main0}
Let $(M^n, g(\t)), \t \in[0, \infty)$ be a complete ancient solution to the Ricci flow with bounded curvature. 
Then for any $p\in M$, we have
\begin{equation}
\lim_{\t \to \infty} \tilV_{(p, 0)}(\t) = \lim_{r \to \infty} I_{(p, 0)}(r).
\end{equation}
\end{thm}

Our Theorem \ref{main0} can be thought of as a general answer to the question quoted above,
as well as a Ricci flow analogue of the following:
for any Riemannian manifold $(M^n, g)$ with non-negative Ricci curvature,
\begin{equation}\label{moti}
\lim_{\t \to\infty} \int_M (4\pi\t)^{-n/2}\exp\({-\frac{d(\cdot, p)^2}{4\t}}\) d\mu
= \lim_{r \to\infty} \frac{\vol B(p, r)}{\omega_n r^n}.
\end{equation}
Here, $\omega_n$ stands for the volume of the unit ball in the Euclidean space $(\mathbb{R}^n, g_E)$.
See \cite[Proposition 12.4]{V2} for the proof.
The right-hand side of (\ref{moti}), denoted by $\nu(g)$, is called the {\it asymptotic volume ratio} of $(M^n, g)$.
The well-definedness of $\nu(g)$ is due to the Bishop-Gromov comparison theorem.

We refer to the limit $\tilNu(g) := \lim_{\t\to\infty}\tilV_{(p, 0)}(\t)$
as the {\it asymptotic reduced volume} of the ancient solution $(M, g(\t)), \t\in [0, \infty)$ to the Ricci flow.
This term was introduced in the author's previous paper \cite{Yo}.
As the notation indicates,
$\tilNu(g)$ is independent of the point $p\in M$ (\cite[Lemma 3.3]{Yo}).
The main theorem of \cite{Yo} was the following gap theorem for ancient solutions to the Ricci flow.
\begin{thm}[{\cite{Yo}}, cf. \cite{Yo2}]\label{gap1}
There exists $\epsilon_n>0$ depending only on $n\ge 2$ such that:
let $(M^n, g(\t)), \t\in [0, \infty)$
be an $n$-dimensional complete ancient solution to the Ricci flow
with Ricci curvature bounded below.
Suppose that its asymptotic reduced volume satisfies $\tilNu(g) > 1-\epsilon_n$.
Then $(M^n, g(\t)), \t \in [0, \infty)$ must be the Gaussian soliton.
\end{thm}

We notice that the assumption on the curvature in our Theorem \ref{main0} is stronger than that in Theorem \ref{gap1}.
This is because there are two different proofs of the monotonicity of the reduced volume in the literature (e.g. \cite{V2});
one of which works only in the bounded curvature setting
where we are able to invoke Shi's estimate.
Shi's gradient estimate serves as a fundamental tool in the study of the Ricci flow,
and our proof of Theorem \ref{main0} presented in this paper also relies on it. 

Instead, we intend to generalize Theorem \ref{main0} in two other directions;
generalizations of the Ricci flow and the reduced volume are discussed in the next section.
See Theorem \ref{main} below for the statement of the main theorem.

The rest of this paper is organized as follows.
In the next section,
we review definitions and results in \cite{P}, \cite{EKNT} and \cite{Ni-Mean}.
This will be done for super Ricci flows.
In Section 3, we restate the main theorem.
The proof of the main theorem is presented in Section 4.

\section{Definitions}
In this section,
we recall the so-called reduced geometry of the Ricci flow introduced by Perelman.
The references are \cite[\S\S 6, 7]{P}, \cite{Ye} and \cite{V2}, etc.
We follow the notation of \cite{V2}.

As mentioned earlier,
we intend to develop the reduced geometry in more general situation.
First,
we recall that a family of Riemannian metrics $(M, g(\t)), \t\in [0, T]$ is called a {\it super Ricci flow} (\cite{MT})
when it satisfies
$$\dertau{}g \le 2\Ric(g(\t)).$$

Apart from the backward Ricci flow,
a basic and important example of a super Ricci flow is
\begin{equation}\label{eq:ex}
g(\t) := (1+2C\t)g_0, \t \in \(0, \frac{1}{|C|-C}\)
\end{equation}
for some fixed Riemannian metric $g_0$ with $\Ric(g_0) \ge Cg_0$ for some constant $C \in \R$.

As was pointed out by the author in \cite{Yo},
it is straightforward to generalize the reduced geometry to the super Ricci flow
if the following assumptions on $2h := \dertau{}g$ are imposed.
\begin{as}\label{as}Letting $H := \tr_{g(\t)} h$, $h$ satisfies
\begin{align}\nonumber\label{evol}
2{\rm div} h(\cdot) &= \langle \nabla H, \cdot\rangle \qquad \text{ (second Bianchi identity)} \text{ and}\\
-\dertau{}H &\ge \varDelta_{g(\t)} H + 2|h|^2 \quad \text{ (heat-like equation)}.
\end{align}
\end{as}

It is clear that
backward Ricci flows ($h=\Ric$) and above examples in (\ref{eq:ex}) ($h=(\frac{1}{C} + 2\t)^{-1}g(\t)$)
are super Ricci flows satisfying Assumption \ref{as}.
The inequality (\ref{evol}), with $\ge$ replaced by $=$,
is nothing but the evolution equation of the scalar curvature $R := \tr_{g(\t)} \Ric$ along the Ricci flow.  
This fact seems to explain a surprising similarity between the two theories
of the Ricci flow and of Riemannian manifolds with non-negative Ricci curvature;
Harnack inequalities, entropy formulae and links to optimal transport theory, and so forth
(e.g. \cite{Mu}, \cite{Ni-Ent}, \cite{To}, \cite{Lo}).

We shall say that a super Ricci flow $(M, g(\t)), \t\in [0, T]$ is {\it $C^1$-controlled}
when we can find a positive function $K(\t)>0$ of $\t$ such that
$$\sup_{M\times [0, \t]} \{|h|+ |\nabla H|^2\} \le K(\t) \text{ for each }\t \in (0, T].$$
According to Shi's gradient estimate (e.g. \cite[Theorem 6.15]{CLN}),
any Ricci flow with bounded curvature is $C^1$-controlled in this sense.

In what follows,
$(M^n, g(\t)), \t\in [0, T]$ or $[0, \infty)$ always denotes a complete $C^1$-controlled super Ricci flow satisfying Assumption \ref{as}
unless otherwise stated.

Next,
we recall the definitions of monotone quantities discovered in \cite{P}, \cite{EKNT} and \cite{Ni-Mean}, respectively.
Fix $p \in M$ and 
let $\phi=\phi(\cdot, \t) \ge0$ be a non-negative locally-Lipschitz function on $M\times [0, T]$.
We assume that the derivatives $\nabla \phi$ and $\dertau{}\phi$ as well as $\phi$ itself satisfy a mild growth condition,
e.g. $\max\{ \phi, |\nabla \phi|^2, |\dertau{}\phi| \} \le \exists C (1+ d_{g(\t)}(\cdot, p)^2)$ for each $\t \in(0, T]$.
 
\begin{defi}[\cite{P}]\label{redvol}
Fix $\bt \in (0, T]$.
We first define the {\it $\L$-length} of a curve $\gamma : \t \mapsto (\gamma(\t), \t) \in M\times [0, \bt]$ by
$$\L(\gamma) := \int_{0}^{\bt} \sqrt\tau \( \|\frac{d\gamma}{d\t}\|^2_{g(\t)} + H(\gamma(\t), \t)\) d\t.$$
Then the {\it reduced distance} between $(p, 0)$ and $(q, \bt) \in M\times (0, T]$ is defined by
$$\l_{(p, 0)}(q, \bt) := \frac{1}{2\sqrt\bt} \inf_{\gamma} \L(\gamma).$$
Here we take the infimum over all curves $\gamma : [0, \bt] \to M$ with $\gamma(0)=p$ and $\gamma(\bt)=q$.
The lower bound for $\dertau{}g$ ensures that the reduced distance is always achieved by a minimal {\it $\L$-geodesic}.

Finally,
letting $K_{(p, 0)}(\cdot, \t) := (4\pi\t)^{-n/2}\exp(-\l_{(p, 0)}(\cdot, \t))$,
we set
$$\tilV_{(p, 0)}^\phi(\t) := \int_M K_{(p, 0)}(\cdot, \t)\phi(\cdot, \t) d\mu_{g(\t)}.$$

When $\phi$ is identically $1$, we drop $\phi$ and
call $\tilV_{(p, 0)}(\t) :=  \int_M K_{(p, 0)}(\cdot, \t) d\mu_{g(\t)}$ the {\it reduced volume} as usual.
The base point $(p, 0)$ is often suppressed.
\end{defi}

\begin{defi}[\cite{EKNT}, \cite{Ni-Mean}]
For any $r>0$,
we let $E_r$ denote the `pseudo heat ball' given by
$$E_r := \{(q, \t) \in M \times (0, T] ; K(q, \t) > r^{-n}\}.$$
Then for a non-negative function $\phi \ge0$ on $M\times [0, T]$, we define
\begin{align*}
I_{(p, 0)}^\phi(r) &:= \frac{1}{r^n} \iint_{E_r} \( |\nabla \log(Kr^n)|^2 + H\log(Kr^n) \) \phi d\mu d\t
\intertext{and}
J_{(p, 0)}^\phi(r) &:= \iint_{\partial E_r} \frac{|\nabla K|^2}{\sqrt{|\nabla K|^2 + |\dertau{}K|^2}}\phi d\tilde{A}
+ \frac{1}{r^n}\iint_{E_r} H\phi d\mu d\t.
\end{align*}
Here $d\tilde{A}$ is the area element induced by the product metric $\tilde{g} := g(\t) + d\t^2$ on $M\times (0, T)$.
\end{defi}

The following proposition stated in \cite{Ni-Mean}
gives the relation between $I_{(p, 0)}^\phi(r)$ and $J_{(p, 0)}^\phi(r)$.
\begin{prop}[\cite{Ni-Mean}]\label{IJ}
For any $r>0$,
$$I_{(p, 0)}^\phi(r) = \frac{n}{r^n} \int_0^r \eta^{n-1} J_{(p, 0)}^\phi(\eta) d\eta.$$
\end{prop}

Let us look at an example.
In our terminology,
Watson's mean value formula for heat equations (\cite{Wa}) can be stated as follows:
\begin{ex}
Let $(\R^n, g(\t) \equiv g_E), \t \in[0, \infty)$ be the Gaussian soliton
and $\phi =\phi(\cdot, \t)$ be a smooth solution to the heat equation $(\dertau{}+\varDelta)\phi =0$.
(Recall that $\t$ is the backward time.)
Take any $p \in \R^n$.
Then $K_{(p, 0)}(\cdot, \t)$ is the heat kernel and
$$I_{(p, 0)}^\phi (r) = \frac{1}{r^ n} \iint_{E_r} |\nabla \log K|^2 \phi d\mu d\t =  \phi(p, 0) \text{ for all } r>0.$$ 
\end{ex}

We now state two theorems asserting that the quantities introduced above are monotone.
\begin{thm}[cf. \cite{P}, \cite{Ye}, \cite{V2}, \cite{Yo}]\label{reduced}
Let $(M^n, g(\t)), \t \in [0, T]$ be a complete $C^1$-controlled super Ricci flow satisfying Assumption \ref{as}.
Suppose that $\phi \ge 0$ satisfies $(\dertau{} + \varDelta_{g(\t)}) \phi \le 0$
in the distributional sense,
namely,
\begin{equation}\label{heat}
\int_M \[ \xi \dertau{\phi} - \langle \nabla \xi, \nabla \phi \rangle \] d\mu_{g(\t)} \le 0
\end{equation}
for any non-negative smooth function $\xi \ge0$ with compact support.
Then for any $p\in M$,
$\tilV_{(p, 0)}^\phi(\t)$ is non-increasing in $\t \in(0, T)$ and $\lim_{\t\to 0+} \tilV_{(p, 0)}^\phi(\t) =\phi(p, 0)$.
\end{thm}

\begin{proof}
The proof of the theorem is identical to the original one for ${d\over d\t} \tilV_{(p, 0)}(\t) \le0$;
the monotonicity of the reduced volume along the Ricci flow (e.g. \cite[Theorem 8.20]{V2}).
See \cite{Yo} for how the original proof is modified for the super Ricci flow satisfying Assumption \ref{as}.
We leave the details to the interested reader.
\end{proof}

\begin{thm}[\cite{EKNT}]\label{thm:EKNT}
Let $(M^n, g(\t)), \t \in [0, T]$ be a complete super Ricci flow satisfying Assumption \ref{as}
with time-derivative $\dertau{}g$ bounded below.
Suppose that $\phi \ge0$ satisfies $(\dertau{}+\varDelta_{g(\t)})\phi \le 0$ in the distributional sense.
For any $p\in M$,
find $r_* >0$ such that $r \in (0, r_*)$ implies that
$E_{r} \subset M\times [0, T-\epsilon)$ for some $\epsilon = \epsilon(r) >0$.
Then
$I_{(p, 0)}^\phi(r)$ is non-increasing in $r \in(0, r_*)$
and $\lim_{r \to 0+} I_{(p, 0)}^\phi(r) = \phi(p, 0)$.
\end{thm}

It was shown by Ni (\cite{Ni-Mean}) that
$J_{(p, 0)}^\phi(r)$ is also non-increasing in $r$
for smooth $\phi \ge0$ and sufficiently small $r>0$ so that $K$ is smooth on $E_r$. 
His point is that
the monotonicity of $I_{(p, 0)}^\phi(r)$ is a consequence of that of $J_{(p, 0)}^\phi(r)$
and Proposition \ref{IJ}.
The following fact is well-known:
$$\text{\it If } \frac{f(r)}{g(r)} \text{ \it is non-increasing in } r >0,
\text{\it then so is } \frac{\int_0^r f(\eta)d\eta}{\int_0^r g(\eta)d\eta}.$$

An example of $\phi \ge0$ as in the above theorems,
other than the constant function $\phi \equiv 1$,
is the function
$$\phi(q, \t) := \max\Bigl\{ 0, \frac{\bar{L}_{(p_*, 0)}(q, \t) -2n\t}{\rho^2} \Bigr\},$$
where $\bar{L}_{(p_*, 0)}(q, \t) := 4\t \l_{(p_*, 0)}(q, \t)$ with $(p_*, 0) \in M\times \{0\}$
and $\rho >0$ is a positive constant.
We deduce this from the inequality $(\dertau{}+\varDelta_{g(\t)}) \bar{L} \le 2n$ (\cite[(7.15)]{P}).
A function which is similar but has compact support was utilized by Ni (\cite{Ni-Matrix})
in order to localize the {\it forward Reduced volume} of Feldman-Ilmanen-Ni (\cite{FIN})
and Perelman's {\it $\mathcal{W}$-functional} (\cite[\S3]{P}).

We refer to \cite{Ec} for earlier works on (local) monotonicity formulae and their applications for mean curvature flow.

We close this section with several facts which will be required later.
They are also utilized in the proofs of Theorems \ref{reduced} and \ref{thm:EKNT}.
The proofs can be found in \cite{Ye}, \cite{V2} etc.

\begin{prop}
If $\dertau{}g \ge -Kg(\t)$ on $M \times [0, \bt]$ for some constant $K \in \R$,
then for any $q \in M$,
$$\l_{(p, 0)}(q, \bt) \ge e^{-K\bt} \frac{d_{g(0)}(p, q)^2}{4\bt} - \frac{nK\bt}{3}.$$ 
\end{prop}

\begin{prop}
For each $\t>0$, we have
$(\dertau{} - \varDelta_{g(\t)} +H)K \le 0$ in the distributional sense.
More precisely,
\begin{equation}\label{comparison}
\int_M \[ -\langle \nabla \l, \nabla \xi \rangle + \( -|\nabla \l|^2+ H + \frac{1}{K}\dertau{K}\)\xi \] d\mu_{g(\t)} \le 0
\end{equation}
holds for any non-negative Lipschitz function $\xi\ge0$ with compact support. 
\end{prop}

\begin{prop}\label{integrable}
Let $(M, g(\t)), \t \in [0, T]$ be a complete $C^1$-controlled super Ricci flow.
Then there exists a positive function $K^*(\t)>0$ of $\t$ such that
$$\max\Bigl\{ |\nabla \l|^2, \Big|\dertau{\l}\Big| \Bigr\} \le \frac{K^*(\t)}{\t}(\l +1) \quad \text{ a.e. on } M$$
for each $\t \in (0, T)$.
In particular, $\int_M |\nabla \l|^2e^{-\l}d\mu_{g(\t)}$ and $\int_M \dertau{\l} e^{-\l}d\mu_{g(\t)}$ make sense.
\end{prop}

\section{Main result}
The main theorem of the present paper will be formulated at the end of this section.
Before that,
let us go back to the question in the introduction.
In order to answer the question,
Ecker et al. (\cite{EKNT}) focused on gradient shrinking Ricci solitons.

\begin{defi}
A triple $(M, g, f)$, with $f$ being a smooth function on $M$,
is called a {\it gradient shrinking Ricci soliton}
if there exists a positive constant $\lambda>0$ such that 
\begin{equation}
\Ric + \Hess f -\frac{1}{2\lambda}g =0.
\end{equation} 
\end{defi}

Given a gradient shrinking Ricci soliton $(M, g, f)$,
we can construct an ancient solution to the Ricci flow (e.g. \cite[Theorem 4.1]{CLN}). 
This is done by expanding and pulling back the metric $g$
by the diffeomorphisms $\phi_\t$ generated by the vector field $-\nabla f$.
The resulting ancient solution $(M, g_0(\t)), \t\in(0, \infty)$ is given by
\begin{equation}\label{soldef}
g_0(\t) := \frac{\t}{\lambda}(\phi_\t)^*g \quad \text{ for } \t\in (0, \infty).
\end{equation}
It was shown by Zhang (\cite{Zh}) that
the gradient vector field $\nabla f$ is complete whenever the metric $g$ is complete.

We define the {\it Gaussian density} (or {\it normalized $f$-volume}) $\Theta(M)$
of gradient Ricci soliton $(M^n, g, f)$ by
$\Theta(M) := \int_M (4\pi\lambda)^{-n/2}e^{-f} d\mu$.
In the sequel,
we tacitly normalize the potential function $f$ by adding a constant so that
$\lambda (|\nabla f|^2 +R) -f \equiv 0$,
with $R$ being the scalar curvature.
This is always possible (e.g. \cite[Lemma 1.15]{V2}).

In this setting, Ecker et al. prove
\begin{prop}[{\cite[Corollary 18]{EKNT}}]\label{EKNT}
Let $(M, g_0(\t)), \t\in(0, \infty)$ be the ancient solution to the Ricci flow
determined by a compact gradient shrinking Ricci soliton $(M, g, f)$.
Fix any $p\in M$.
Then for any $\t>0$ and $r>0$,
$$\tilV_{(p, 0)}(\t) = I_{(p, 0)}(r) = \Theta(M).$$
\end{prop}

Note that in the situation of Proposition \ref{EKNT},
$(M, g_0(\t))$ shrinks to a point as $\t \to 0+$ and
what they are dealing with is the reduced distance form the {\it singular} point.
Hence the results like
\begin{itemize}
\item {\it $\Theta(M) \le 1$ and
\item $\Theta(M) =1$ $\iff$ $(M^n, g)$ is isometric to $(\R^n, g_E)$}
\end{itemize}
does not follow from Proposition \ref{EKNT} (immediately, at least).
It was \cite{CN} and \cite{Yo}
where such comparison geometric results were established under respective curvature conditions.
(See Proposition \ref{redvol=fvol} below.)
This is how our Theorem \ref{main0} differs from Proposition \ref{EKNT}.
(We should mention \cite{Na} and \cite{E}
where the reduced volume based at a singular point was studied under some singularity assumptions.)

Before stating our main theorem,
let us prove the following
\begin{prop}\label{prop:stat}
Let $(M^n, g)$ be a complete Riemannian manifold of non-negative Ricci curvature regarded as a static super Ricci flow,
i.e., $g(\t)\equiv g$.
Then for any $p \in M$,
$$\lim_{r\to\infty} I_{(p, 0)}(r) = \nu(g).$$
Here $\nu(g)$ denotes the asymptotic volume ratio of $(M^n, g)$ as before.  
\end{prop}

\begin{proof}
As was shown in \cite[Lemma 9]{EKNT}, with $\psi := \log K$,
\begin{equation}\label{altform}
I_{(p, 0)}(r) = \frac{1}{r^n} \iint_{E_r} \[ |\nabla \psi|^2 - \dertau{\psi} \] d\mu d\t \text{ for all } r>0.
\end{equation}
Using (\ref{altform}) and that $\l_{(p, 0)}(\cdot, \t) = \frac{1}{4\t}d(\cdot, p)^2$,
we know
$$I_{(p, 0)}(r)
= \frac{1}{r^n} \iint_{E_r} \frac{n}{2\t} d\mu d\t
= \frac{1}{r^n} \int_0^{\frac{r^2}{4\pi}} \frac{n}{2\t} \vol B\(p, \sqrt{2n\t \log\frac{r^2}{4\pi\t}}\) d\t.$$
The rest of the proof is a simple calculation.
\end{proof}

Now we are in a position to state the main result of this paper.
\begin{thm}\label{main}
Let $(M^n, g(\t)), \t \in[0, \infty)$ be a super Ricci flow which is complete, ancient and $C^1$-controlled. 
Then for any $p\in M$ and a non-negative locally-Lipschitz function $\phi \ge0$
with $(\dertau{} + \varDelta_{g(\t)}) \phi \le 0$ in the distributional sense,
we have
\begin{equation}\label{eq:main}
\lim_{\t \to \infty} \tilV^\phi_{(p, 0)}(\t) = \lim_{r \to \infty} I_{(p, 0)}^\phi(r).
\end{equation}
\end{thm}
Clearly, Theorem \ref{main} contains Theorem \ref{main0} as a special case.

\section{Proof}
\begin{proof}[Proof of Theorem \ref{main}]
In order to establish (\ref{eq:main}),
we follow the same line as in the proof of Theorem \ref{thm:EKNT} given in \cite{EKNT}.

First of all, we fix small $\epsilon >0$.
Take a $C^\infty$-function $\eta :(-\infty, \infty) \to [0, \infty)$
such that the support is contained in $[0, \epsilon]$ and $\int_0^\infty \eta(y)dy =1$.
We define $\zeta(x) := \int_0^x \eta(y)dy$ and $Z(x) := \int_0^x \zeta(y)dy$.
Notice that $\eta$ is a smooth approximation of the Delta function, and hence,
$\zeta$ and $Z$ approach to the Heviside function $\chi$
and the function $x \mapsto [x]_+ :=\max\{x, 0\}$ as $\epsilon \to 0+$, respectively.
More precisely, for any $x \in (-\infty, \infty)$, we have
$$\chi(x-\epsilon) \le \zeta(x) \le \chi(x)$$
and hence
$$[x-\epsilon]_+ \le Z(x) \le [x]_+.$$

Set
$$ Q(s, r) := \iint_{M\times [s, \infty)} \( |\nabla \log K|^2 \zeta(\log(Kr^n)) +HZ(\log(Kr^n)) \) \phi d\mu d\t$$
for all $s \ge0$ and $r> 0$.
We also put $\I(s, r) := Q(s, r)/r^n$ and $\I(r) := \I(0, r)$.
Then we know
\begin{equation}\label{error}
e^{-\epsilon} I_{(p, 0)}^\phi(e^{-\epsilon/n}r) \le \I(r) \le I_{(p, 0)}^\phi(r) \ \ \text{ for all } r>0.
\end{equation}
Here we used the fact (\cite{Ch}, \cite[Proposition A.2]{Yo}) that:
{\it For any complete ancient super Ricci flow} $(M, g(\t))$ {\it satisfying} (\ref{evol}), $H$ {\it is non-negative on} $M\times [0, \infty)$.

Next, we have
\begin{align*}
&{d\over dr} \I(s, r)\\
=& \frac{n}{r^{n+1}} \(\frac{r}{n} {d\over dr}Q(s, r) - Q(s, r)\)\\
=& \frac{n}{r^{n+1}} \iint_{M\times [s, \infty)} \( |\nabla \log K|^2 \zeta^\prime + HZ^\prime - |\nabla \log K|^2\zeta -HZ \) \phi d\mu d\t.
\end{align*}

We use the well-known formula $\dertau{} d\mu_{g(\t)} =Hd\mu_{g(\t)}$ to get
$${d\over d\t} \int_M Z(\log(Kr^n))\phi d\mu_{g(\t)}
= \int_M \[ \(Z^\prime {1\over K}\dertau{K} +ZH\) \phi + Z \dertau{\phi} \] d\mu_{g(\t)}.$$
Since $\int_{M\times \{s^\prime\}} Z\phi d\mu =0$
for all $s^\prime>>1$ with $M \times [0, s^\prime) \supset E_r$,
we obtain
$$- \int_{M\times \{s\}} Z(\log(Kr^n))\phi d\mu
= \iint_{M\times [s, \infty)} \[ \(Z^\prime {1\over K}\dertau{K} +ZH\) \phi + Z \dertau{\phi} \] d\mu d\t.$$
Thus
\begin{equation}\label{eq:concl}
{d\over dr} \I(s, r)
= \frac{n}{r^{n+1}} \iint_{M\times [s, \infty)} \mathcal{A}d\mu d\t
+ \int_{M\times \{s\}} Z(\log(Kr^n)) \phi d\mu,
\end{equation}
where we let
$$ \mathcal{A} := \( |\nabla \log K|^2 \zeta^\prime + HZ^\prime - |\nabla \log K|^2\zeta + Z^\prime {1\over K}\dertau{K} \) \phi + Z\dertau{\phi}.$$

We also let
\begin{align*}
\mathcal{A}_* :=& -\langle\nabla \l, \nabla (\zeta\phi)\rangle + \(-|\nabla \l |^2 +H + \frac{1}{K}\dertau{K}\)(\zeta\phi)\\
&\qquad + \(Z\dertau{\phi} -\langle \nabla Z, \nabla \phi \rangle\)
\end{align*}
to observe that
\begin{equation}\label{key}
\int_M \mathcal{A} d\mu_{g(\t)} = \int_M \mathcal{A}_* d\mu_{g(\t)} \le 0 \text{ for all } \t>0.
\end{equation}
To show (\ref{key}), we have used
\begin{align*}
\int_M |\nabla\log K|^2\zeta^{\prime} \phi d\mu
&= \int_M \langle \nabla\log K, \nabla\zeta \rangle \phi d\mu\\
&= \int_M \[ \langle \nabla\log K, \nabla(\zeta\phi) \rangle - \zeta \langle \nabla \log K, \nabla\phi \rangle \] d\mu\\
&= \int_M \[ \langle \nabla\log K, \nabla(\zeta\phi) \rangle - \langle \nabla Z, \nabla\phi \rangle \] d\mu. 
\end{align*}
The last inequality in (\ref{key}) follows from (\ref{heat}) and (\ref{comparison}).

Integrating (\ref{eq:concl}) for $r_1 \le r \le r_2$ yields
\begin{align*}
&\I(s, r_2)- \I(s, r_1)\\
=& \int_{r_1}^{r_2} \frac{n}{r^{n+1}}dr \int_s^\infty d\t \int_M \mathcal{A} d\mu_{g(\t)}  + E(s; r_1, r_2)\\
=& \int_s^\infty d\t \int_{r_1}^{r_2} \frac{n}{r^{n+1}}dr \int_M \mathcal{A}_* d\mu_{g(\t)} + E(s; r_1, r_2),
\end{align*}
where $E(s; r_1, r_2)$ denotes the extra term:
$$E(s; r_1, r_2) := \int_{r_1}^{r_2} dr \int_{M\times \{s\}} Z(\log(Kr^n)) \phi d\mu.$$
We have applied Fubini's theorem.
This can be done freely due to (\ref{key}). 

Now we see that $E(s; r_1, r_2) \to 0 \text{ as } s \to 0+$.
Indeed,
\begin{align*}
&\limsup_{s\to 0+} E(s; r_1, r_2)\\
\le& \limsup_{s\to 0+}\, (r_2-r_1)
\(n\log\frac{r_2}{\sqrt{4\pi s}}\) \vol_{g(s)} \Bigl\{ \l(\cdot, s) < n\log \frac{r_2}{\sqrt{4\pi s}} \Bigr\}\\
\le& \limsup_{s\to 0+}\, (r_2-r_1)
\(n\log\frac{r_2}{\sqrt{4\pi s}}\) \omega_n \(4ns \log\frac{r_2}{\sqrt{4\pi s}}\)^{n/2}\\
=&0.
\end{align*}
With this observation, 
letting $s\to 0+$ yields
\begin{equation}\label{eq:mono}
 \I(r_2)- \I(r_1)= \int_0^\infty d\t  \int_{r_1}^{r_2} \frac{n}{r^{n+1}}dr \int_M  \mathcal{A}_* d\mu_{g(\t)} \le 0,
\end{equation}
which implies the monotonicity of $\I(r)$ in $r>0$.
(We can now let $\epsilon \to 0+$ to see that $I_{(p, 0)}^\phi(r)$ is non-increasing in $r>0$.)

Here, for any positive $K>0$, we set
$$K_\eta := \int_0^\infty \frac{n}{r^{n+1}} \eta(\log(Kr^n))dr,$$
and $K_\zeta$ and $K_Z$ are also defined similarly.
It is easy to check by using the integration by parts that
$$K_\eta = K_\zeta = K_Z = e^{\delta(\eta)}K, \text{where } e^{\delta(\eta)} := \int_0^\infty \eta(y)e^{-y}dy.$$
Notice that $\delta(\eta) \le 0$ and $\delta(\eta) \to 0$ as $\epsilon \to 0+$.

Then for any $\t>0$,
\begin{align*}
&\int_0^\infty \frac{n}{r^{n+1}}dr \int_M \mathcal{A} d\mu_{g(\t)}\\
=& \int_M \[\( |\nabla \l|^2 K_\eta + H K_\zeta -|\nabla \l|^2 K_\zeta + \frac{K_\zeta}{K}\dertau{K} \)\phi + K_Z \dertau{\phi} \]d\mu_{g(\t)}\\
=& e^{\delta(\eta)} \int_M \[ \( |\nabla \l|^2 K + H K  -|\nabla \l|^2 K + \dertau{K} \)\phi + K \dertau{\phi} \] d\mu_{g(\t)}\\
=& e^{\delta(\eta)} \int_M \[\( H K + \dertau{K} \)\phi + K \dertau{\phi}\] d\mu_{g(\t)}.
\end{align*}
We are implicitly using the integrability of each term (Proposition \ref{integrable}) to derive the equations above.

By letting $r_2 \to \infty$ and $r_1 \to 0$ in (\ref{eq:mono}), we then get
\begin{align*}
&\lim_{r\to\infty}\I(r) - \lim_{r\to0}\I(r)\\
=&\int_0^\infty d\t \int_0^\infty \frac{n}{r^{n+1}}dr \int_M \mathcal{A} d\mu_{g(\t)}\\
=& e^{\delta(\eta)} \int_0^\infty d\t \int_M \[\( H K + \dertau{K} \)\phi + K \dertau{\phi}\] d\mu_{g(\t)}.
\end{align*}

Finally, we take $\epsilon \to 0+$ and use (\ref{error}) to conclude that
\begin{align}\label{1}
\lim_{r\to\infty} I_{(p, 0)}^\phi(r) - \phi(p, 0)
= \int_0^\infty d\t \int_M \[\( H K + \dertau{K} \)\phi + K \dertau{\phi}\] d\mu_{g(\t)}.
\end{align}

On the other hand, we know that
\begin{align}\label{2}
\lim_{\t \to \infty}\tilV^\phi_{(p, 0)}(\t) - \phi(p, 0)
= \int_0^\infty d\t \int_M  \[\(\dertau{K} +KH\)\phi + K\dertau{\phi}\]d\mu_{g(\t)}.
\end{align}
Combining (\ref{1}) and (\ref{2}) completes the proof of Theorem \ref{main}.
\end{proof}

Let us reconsider shrinking Ricci solitons.
We let $g_S(\t) := g_0(\t +s), \t\in [0, \infty)$ for some fixed $s>0$,
where $g_0(\t)$ is the ancient solution defined in (\ref{soldef}).
Notice that $\t=0$ is no longer the singular time for $(M, g_S(\t))$. 
\begin{prop}[{\cite[Proposition 5.1]{Yo}}]\label{redvol=fvol}
Let $(M, g, f)$ be a complete gradient shrinking Ricci soliton with Ricci curvature bounded below.
Then
$$\Theta(M) = \tilNu(g_S),$$
that is, the Gaussian density $=$ the asymptotic reduced volume.
 \end{prop}

Combined with the above proposition,
Theorem \ref{main0} implies the following corollary (compare with Proposition \ref{EKNT}).
\begin{cor}
Let $(M, g_S(\t)), \t\in [0, \infty)$ be the ancient solution to the Ricci flow
determined by a complete gradient shrinking Ricci soliton $(M, g, f)$ with bounded curvature.
Then for any $p \in M$,
$$\lim_{\t\to \infty}\tilV_{(p, 0)}(\t) = \lim_{r\to\infty} I_{(p, 0)}(r) = \Theta(M).$$ 
\end{cor}

We conclude this paper with a few remarks.
\begin{rem}
\begin{enumerate}
\item The author wonders whether Theorem \ref{main0} still holds under the assumption of Theorem \ref{gap1}.
He believes that a more understanding of Perelman's reduced geometry
from a geometric viewpoint is required to attack the problem.
The reduced volume $\tilV(\t)$ and $I(r)$ are well-defined
as long as $\dertau{}g$ is bounded from below (see \cite{Ye}, \cite{EKNT}).
\item The reader is referred to Ni's paper \cite{Ni-Ent}
where he raises an interesting question closely related to our Theorem \ref{main0}.
\end{enumerate}
\end{rem}

\section*{Acknowledgements}
The author would like to express his gratitude to his adviser Takao Yamaguchi for his encouragement.
This work was supported in part by Research Fellowships of the Japan Society for the Promotion of Science for Young Scientists.

\end{document}